\documentclass[a4paper,12pt,english]{amsart}

\usepackage[english]{babel}
\usepackage[latin1]{inputenc}
\usepackage[T1]{fontenc}

\usepackage{graphicx}
\usepackage{amssymb}
\usepackage{verbatim,enumerate}

\newtheorem{thm}{Theorem}[section]

\newtheorem{prop}[thm]{Proposition}
\newtheorem{cor}[thm]{Corollary}
\newtheorem{defi}[thm]{Definition\rm}

\newcommand{\C}{\mathbb{C}}
\newcommand{\CP}{\mathbb{C}P}
\newcommand{\RP}{\mathbb{R}P}
\newcommand{\F}{\mathbb{F}}

\newcommand{\R}{\mathbb{R}}
\newcommand{\Z}{\mathbb{Z}}
\newcommand{\PC}{\widetilde{\mathbb C P^2_6}}
\newcommand{\PR}{\widetilde{\mathbb R P^2_6}}
\newcommand{\PP}{\mathbb C P^2_6}
\newcommand{\PPR}{\mathbb R P^2_6}
\newcommand{\x}{\underline{x}}

\begin{document}

\title{Behavior of Welschinger invariants under Morse simplifications}

\author{Erwan Brugall\'e}
\address{Université Pierre et Marie Curie,  Paris 6, 4 Place Jussieu, 75~005 Paris, France}

\email{brugalle@math.jussieu.fr}

\author{Nicolas Puignau}
\address{Universidade Federal do Rio de Janeiro, Ilha do Fundão, 21941-909 Rio de Janeiro, Brasil}

\email{puignau@im.ufrj.br}


\maketitle

\begin{abstract}
We relate Welschinger invariants of a rational real symplectic 4-manifold
before and after a Morse simplification (i.e deletion of a
sphere or a handle of the real part of the surface). This relation is 
a consequence of a real version of Abramovich-Bertram formula 
 which computes Gromov-Witten
invariants  by means of enumeration of
$J$-holomorphic curves with a non-generic almost complex structure $J$.  
In addition, we give some qualitative consequences of our study, for
example the vanishing of  Welschinger invariants in some cases.
\end{abstract}

\section{Introduction}
On a rational symplectic 4-manifold $(X,\omega)$, 
genus 0 Gromov-Witten invariants
can be computed by enumerating \textit{irreducible} $J$-holomorphic
rational curves on $X$, 
realizing a fixed  homology class $d\in H_2(X,\mathbb Z)$, 
and passing through a  configuration of $c_1(X)d -1$
points, where $J$ is a generic almost complex structure on $X$ tamed
 by $\omega$ (\cite{McS}).
Now suppose that $J$ is midly non-generic, i.e. $X$ contains a unique
irreducible 
$J$-holomorphic curve $E$ with $E^2< -1$, and moreover  $E$ is
a smooth rational curve with $E^2=-2$. In this situation, one can still compute Gromov-Witten
invariants of $(X,\omega)$ by 
enumerating $J$-holomorphic curves on $(X,\omega)$, but now 
also taking
into account  \textit{reducible} curves with some components mapped
isomorphically to $E$. Abramovich and Bertram first proved this
when $(X,\omega,J)$ is the second
Hirzebruch ruled surface (\cite{AB}), Vakil extended later this proof to
the case of any 
weak Del Pezzo surface (\cite{V}), and eventually
 Ionel and Parker symplectic sum
formula (\cite{IP}) provides a proof in the general case.

Results of this note are based on real versions of this
Abramovich-Bertram type formula. 
A real structure $c:X\to X$ on a rational symplectic 4-manifold
$(X,\omega)$
is an involution such that $c^*\omega=-\omega$. The set $\R X=Fix(c)$ is
called the
real locus of $X$. 
Welschinger
 invariants provide real analogues
 of Gromov-Witten invariants in genus 0 for real rational symplectic 4-manifolds
 (\cite{W}). 

\vspace{1ex}
Suppose that $(X, \omega,c)$ contains a real smooth rational symplectic curve $E$
with $E^2=-2$, and let $(X^\#,\omega^\#)$ be the symplectic sum of
$(X,\omega)$ with $S^2\times S^2$ along $E$, where $E$ realizes the
diagonal class in $H_2(S^2\times S^2,\Z)$. 
There exist 
two real structures $c_+$ and $c_-$
on $S^2\times S^2$
for which $E$ is real, which give rise to two different real
structures $c_+^\#$ and $c_-^\#$ on $(X^\#,\omega^\#)$
satisfying 
 (with the
 convention that 
 $\chi(\emptyset)=0$)
$$\chi(\R X_+^\#)=\chi(\R X)=\chi(\R X_-^\#)-2.$$
One may interpret this construction as follows: 
blow-down the real
(-2)-curve $E$ to a nodal real 4-manifold, and  
smooth the node in two different ways.

The real symplectic manifold $(X^\#,\omega^\#,c_+^\#)$ is in fact a
 deformation of $(X,\omega,c)$ and in this case one can
immediatly
extract  a real 
version of Abramovich-Bertram formula from the complex one without
decomposing $(X^\#,\omega^\#,c_+^\#)$ into a symplectic sum, 
as it has already been
noticed by several people (\cite{Br2}, \cite{Br}, \cite{Kha}, \cite{RasSal}).
This is not true for $(X^\#,\omega^\#,c_-^\#)$, and
one of the main results of this note is 
a real version of Abramovich-Bertram formula also in this case.
These two different 
real versions of  Abramovich-Bertram formula
allows one to compare 
Welschinger invariants of
$(X^\#,\omega^\#,c_\pm^\#)$.
This can be thought as a generalization of the 
invariant $\theta$ introduced by Welschinger in \cite{W}, and has 
several consequences (e.g. vanishing results) concerning Welschinger
invariants.

\vspace{1ex}
Detailed proofs of the statements announced in this note will appear
in \cite{BP}.

\vspace{1ex}
\textbf{Acknowledgement: }This final formulation of our results
wouldn't have been possible without the patient explanations of many of
our colleagues. We are in particular
indebted to Simone Diverio, Ilia Itenberg,
Viatcheslav Kharlamov, Frédéric Mangolte, Christian
Peskine, Brett Parker, Patrick Popescu, and Jean-Yves Welschinger.

Both authors were supported by the
Brazilian-French Network in Mathematics. E.B. is also
partially supported by the ANR-09-BLAN-0039-01 and
ANR-09-JCJC-0097-01.

\section{Welschinger invariants}\label{def}

Let $(X,\omega,c)$ be a real rational symplectic 4-manifold, and let
$J$ be an almost complex structure on $X$ tamed by $\omega$ which
is $J$-antiholomorphic.
Recall that the \textit{mass} $m(C)$ of a real rational  $J$-holomorphic
curve $C$ in  $(X,\omega,c)$ is the number of
solitary real nodes of $\mathbb RC$ in $\mathbb R X$ (i.e. nodes
locally given over $\mathbb R$ by the equation $x^2+y^2=0$). Let us
fix a homology class $d$ in $H_2(X,\mathbb Z)$, an integer $0\le r\le
c_1(X)d-1$, a connected component $S$ of $\R X$, and a 
real configuration $\x$ of $c_1(X)d-1$ points in $X$ containing
exactly $r$ points in $S$
 and $\frac{c_1(X)d-1-r}{2}$ pairs of
complex conjugated points.  
When $J$ is generic, Welschinger proved in \cite{W} that the
number of irreducible real rational $J$-holomorphic curves $C$, counted with
multiplicity $(-1)^{m(C)}$, incident to $\x$ and realizing the
class $d$ is finite and 
depends only on $d$ and $r$.
This number is a \textit{Welschinger invariant} of $(X,\omega,c)$, and we
denote it by 
$W_{\R X,S} (d,r)$.
We omit the reference to $S$ when $S=\R X$, or to $r$ 
when $r=c_1(X)d-1$. 

Suppose now that $J$ is mildly non-generic as above, in particular the
$(-2)$-curve $E$ is real.
Counting real rational $J$-holomorphic  curves in 
$X$ with multiplicity  $(-1)^{m(C)}$ does not give a number
depending only on $d$ and $r$,
 since  $J$ is non-generic (\cite{W}, \cite{IKS3}).

\begin{defi}
Let $C$ be a nodal real rational $J$-holomorphic  curve in 
$X$ intersecting the (-2)-curve $E$ transversally. 
We denote respectively by $\alpha$ and $\beta$ the  number of real and
pairs of complex conjugated intersection points in $C\cap E$. For any
integer $k\ge 0$, we define the two $k$th multiplicities  of $C$ as follows:
$$\mu^+_k(C)=(-1)^{m(C)}\sum_{k=\alpha_k+2\beta_k}{\alpha \choose \alpha_k}{\beta\choose \beta_k}$$
and
$$\mu^-_k(C)=\left\{ \begin{array}{ll}
         (-1)^{m(C)+\beta}2^\beta & \mbox{if $\alpha=0$ and $k=\beta$};\\
         0 & \mbox{otherwise}.\end{array} \right.$$

\end{defi}
 As above choose $d\in H_2(X,\mathbb Z)$,
 an integer $0\le r\le c_1(X)d-1$, a connected component $S$ of
 $\mathbb RX \setminus \mathbb RE$,
 and a generic real configuration $\x$ of $c_1(X)d-1$
  points in $X$ containing exactly $r$ points in $S$ and
$\frac{c_1(X)d-1-r}{2}$ pairs of complex conjugated points. 
For each integer $k\ge 0$, we denote by $\mathcal R_k(d,\omega)$
the set of all irreducible rational real $J$-holomorphic curves in $X$
passing through all points in $\x$ and realizing the class
$d-kE$. 
The set
$\mathcal R_k(d,\omega)$ is finite, 
and any curve in $\mathcal R_k(d,\omega)$ is nodal and intersects
$E$ transversally. 
Moreover $\mathcal R_k(d,\omega)$ is non-empty only for finitely
many values of $k$.
We define 
the two following numbers:
$$W^\pm_{\R X,S}(d,r)=\sum_{k\ge 0}
\  \sum_{C\in\mathcal R_k(d,\omega)}
\mu^{\pm}_k(C).$$
Let $(X^\#,\omega^\# ,c^\#)$ 
be as above with $c^\#=c_\pm^\#$,
 and let $S^\#$ be the
component of $\R X^\#$ 
containing the deformation of $S$. 
Note that
the homology groups $H_2(X,\mathbb Z)$ and $H_2(X^\#,\mathbb
Z)$ are canonically identified (\cite{IP}).

\begin{thm}\label{mainthm}
Under the above hypotheses,  one has:
\begin{enumerate}[(i)]
\item if $\chi(\R
X^\#)=\chi(\R X)$, then
$$W_{\R X^\#,S^\#}(d,r)=W^+_{\R X,S}(d,r);$$
\item if  $\chi(\R
X^\#)=\chi(\R X)+2$, then
$$W_{\R X^\#,S^\#}(d,r)=W^-_{\R X,S}(d,r).$$
\end{enumerate}
\end{thm}

As an immediate consequence of Theorem \ref{mainthm}, the numbers
$W^\pm_{\R X,S}(d,r)$ 
depend 
only on $d$ and $r$. 
As mentioned in the introduction, part $(i)$ in Theorem \ref{mainthm}
is an immediate consequence of Abramovich-Bertram formula and
was known before (\cite{Br2},
\cite{Br}, \cite{Kha}, \cite{RasSal}).

\section{Applications}

Here we announce some 
consequences of
Theorem
\ref{mainthm}, in particular when $X$ is $\PP$, the complex projective
plane $\CP^2$ blown up in 6 points. 

\subsection{Computation for degree 6 curves with 6 fixed nodes} 
Let us also denote by $\PC$ the
projective plane $\C P^2$  blown up at 6
points lying on a smooth conic $E$.
Here we enumerate  real rational curves 
realizing twice the anti-canonical class $\delta=2c_1(\PP)^{\vee}$ in
$\PC$ and $\PP$. 

Given a real structure on $\PC$, we denote by $\PR$ its real
part. Note that $\PR$ is not necessarily $\RP^2$ blown up in 6 real
points lying on a conic.
Given a generic
configuration $\x$ of $c_1(\PP)\delta-1=5$ real points in $\PC$, 
we set $n_{\chi(\PR)}^{\pm}(\delta-kE):=\sum_{C\in\mathcal
  R_k(\delta,\omega)}\mu^{\pm}_k(C)$.  

\begin{prop}\label{actual comp}
For any choice of $S$, there exists a configuration of 5 real
points in $\PC$ such that:\begin{center}
\begin{tabular}{ c|c| c|c|c| c|c|c|c}
   & $n^{+}_{-5}$& $n^{-}_{-5}$& $n^{+}_{-3}$& $n^{-}_{-3}$& $n^{+}_{-1}$&$n^{-}_{-1}$&
  $n^{+}_{1}$&   $n^{-}_{1}$
\\\hline
  $\delta$   & 522 & 522 & 236 & 236 & 78 & 78 & 0 &0
\\\hline
 $\delta-E$ & 236 & 0 & 140 & 0 & 76 & 0 & 36 & 0
\\\hline
  $\delta-2E$ & 1&0&1&0&1&0&1&0
\end{tabular}
\vspace{1ex}

\end{center}
\end{prop}

\begin{cor}
The surface $\PP$ has the following Welschinger invariants: 
\begin{center}
\begin{tabular}{c| c|c|c|c|c}
 $\chi(\PPR) $&    $-5$    &   $-3$   &  
$-1$
   &  $1$     & $3$   
\\
\hline

$W_{\PPR,S}(\delta)$  &  1000 & 522& 236  &78 &0 
\end{tabular}
\end{center}
\end{cor}

The value $W_{\PPR}(\delta)$  when $\chi(\PPR)=-5$
has been first
    computed by the first author (\cite{Br2}, \cite{Br}). The numbers
  $W_{\PPR}(\delta)$  when $\chi(\PPR)=-3,-1,1$, as well as 
 $W_{\R P^2\sqcup S^2,\R P^2}(\delta)$
have
              been first computed by Itenberg, Kharlamov and Shustin
              (\cite{IKS11}). 
The vanishing of $ W_{\R P^2\sqcup S^2, S}(\delta)$ 
is actually a general fact. 

\begin{prop}\label{null} If $(X,\omega,c)$ is a real symplectic
    4-manifold 
with  disconnect real part, then for any  $d\in H_2(X, \Z)$, any $r\ge 2$,
and any choice of $S$, one has
$$W_{\R X, S}(d,r)=0.$$
\end{prop}


\subsection{Behavior of purely real Welschinger invariants with
  respect to Euler characteristic}

Given a real toric Del Pezzo surface $X$
equipped with its 
tautological real toric structure and
a class $d\in H_2(X,\Z)$, one has
(\cite{IKS04})
$$W_{\R X}(d)\ge W_{\R X}(d,c_1(X)d-3).$$
Theorem \ref{mainthm} provides a natural generalization of this
formula in the particular cases when $X$ is  $S^2\times S^2$ or
$\PP$. 

\begin{thm}\label{decr}
Let $(X_1,\omega_1)$ and $(X_2,\omega_2)$ 
 be two symplectic 4-manifolds deformation equivalent to either 
$\C P^1\times \C P^1$ or
$\PP$ equipped with their standard symplectic form. 
Choose a real structure $c_1$ on $X_1$, and a real structure $c_2$ on
$X_2$. 
 Then
for any $d\in H_2(X,\Z)$, one has
$$W_{\R X_1,S_1}(d)\ge  W_{\R X_2,S_2}(d)\quad \text{if}\quad
\chi(\R X_1)\le \chi(\R X_2).$$
\end{thm}

Note that Theorem \ref{decr} does not generalize immediately to any
symplectic 4-manifold. Indeed, according to \cite{ABL} one has
$W_{\R P^2}(9, 2)< W_{\R P^2}(9, 0)$, i.e. Theorem \ref{decr} does not
hold in the case of $\C P^2$ blown up in 26 points.

\subsection{Modified Welschinger invariants}
In the case when $\R X$ is not connected, one may slightly modify the
definition of Welschinger invariants given in section \ref{def}. 
Namely, given $S$ a connected component of $\R X$,
the modified mass of a real rational
curve $C$ is defined as the number of solitary real nodes of $C$ lying
in $S$. 
Counting real curves with this sign produces a new invariant, denoted
by $\widetilde W_{\R X,S}$.

Our method also allows us to compute these invariants in the case of
$\PP$. 
In particular we have the following two propositions.
\begin{prop} 
$\widetilde W_{\RP^2\sqcup S^2,\R P^2}(\delta)=160 \quad \textrm{and}\quad 
\widetilde W_{\RP^2\sqcup S^2,S^2}(\delta)=96.$
\end{prop}
The value of $\widetilde{W}_{\RP^2\sqcup S^2,\R P^2}(\delta)$ has
              been first computed by Itenberg, Kharlamov and Shustin
              (\cite{IKS11}).

\begin{prop}
For any  class $d\in H_2(\PP,\Z)$, we have
$$\widetilde W_{\RP^2\sqcup S^2,\R P^2}(d)\ge
\widetilde W_{\RP^2\sqcup S^2,S^2}(d)\ge 0. $$
\end{prop} 

The positivity of $\widetilde W_{\RP^2\sqcup S^2,\R P^2}(d)$ whenever
$d$ contains a real algebraic curve has first been established in \cite{IKS11}. 

\subsection{Relation to tropical Welschinger invariants of $\mathbb
  F_2$}\label{sec:torp F2}
We end this note  relating some tropical
Welschinger invariants of $\F_2$ to genuine Welschinger invariants of
the quadric ellipsoid $Q$. The only real homology classes  of $Q$ are
multiple of the hyperplane section $h$.
We say that a tropical curve in $\mathbb R^2$ is of class $aB+bF$ in
$\mathbb T\mathbb F_2$ if
its Newton polygon has vertices $(0,0)$, $(0,a)$, $(a,b)$, and
$(2a+b,0)$. 
We denote by $W_{\mathbb T\mathbb  F_2}(dB)$ the irreducible tropical
Welschinger invariant of $\mathbb T\mathbb
F_2$ for curves of class $dB$ (\cite{IKS09}).

\begin{prop}\label{class-trop}
For any positive integer $d$ 
$$W_{Q}(dh)=W_{\mathbb T
 \mathbb F_2}(dB). $$
\end{prop}

\end{document}